\documentstyle[12pt, righttag]
{amsart}
\input amssym.def
\input amssym

\newcommand{\wt}{\mbox{\normalshape wt}}
\newcommand{\spa}{\mbox{\normalshape span}}
\newcommand{\Res}{\mbox{\normalshape Res}}

\newcommand{\Ind}{\mbox{\normalshape Ind}}
\newcommand{\Hom}{\mbox{\normalshape Hom}}

\def \Z{\Bbb Z}

\def \C{\Bbb C}

\def \N{\Bbb N}

\def \<{\langle} 

\def \O{\Omega}

\def \1t{\frac{1}{T}}
\def \>{\rangle}

\def \l{\lambda }

\def \be{\begin{equation}\label}
\def \bee{\begin{equation}\label}
\def \ee{\end{equation}}
\def \qed{\mbox{ $\square$}}
\def \pf {\noindent {\bf Proof:} \,}
\def \bl{\begin{lem}\label}
\def \el{\end{lem}}
\def \ba{\begin{array}}
\def \ea{\end{array}}
\def \bt{\begin{thm}\label}
\def \et{\end{thm}}

\def \br{\begin{rem}\label}
\def \er{\end{rem}}
\def \ed{\end{de}}
\def \bp{\begin{prop}\label}
\def \ep{\end{prop}}

\newtheorem{thm}{Theorem}[section]
\newtheorem{prop}[thm]{Proposition}

\newtheorem{lem}[thm]{Lemma}
\newtheorem{rem}[thm]{Remark}

\numberwithin{equation}{section}

\begin{document}
\title[Vertex operator algebras and associative algebras]{vertex operator algebras and associative algebras}
\author{Chongying Dong, Haisheng Li and Geoffrey Mason}
\address{Department of Mathematics,\,University of California,\,Santa Cruz (C.D., G.M.)\\ Department of Mathematical Sciences,
Rutgers University, Camden, NJ 08102 (H.L.)}
    \email{dong@@cats.ucsc.edu\\
    hli@@crab.rutgers.edu\\
    gem@@cats.ucsc.edu}
\thanks{C.D. is partially supported by NSF grant DMS-9303374 and a
research grant from the Committee on Research, UC Santa Cruz.}
\thanks{H.L. is partially supported by NSF grant DMS-9616630.}
\thanks{G.M. is partially supported by NSF grant
DMS-9401272 and a research grant from the Committee on Research, UC
Santa Cruz.}
\bibliographystyle{alpha} 
	\maketitle

\begin{abstract} Let $V$ be a vertex operator algebra.
We construct a sequence of associative
algebras $A_n(V)$ ($n=0,1,2,...$) such that $A_{n}(V)$ is a quotient of
$A_{n+1}(V)$ and a pair of functors between
the category of $A_n(V)$-modules which are not $A_{n-1}(V)$-modules
and the category of admissible $V$-modules. These
functors exhibit a bijection between the simple modules in each
category. We also show that $V$ is rational if and only if all 
$A_n(V)$ are finite-dimensional semisimple algebras.
\end{abstract}

\section{Introduction}

For a vertex operator algebra $V$ Zhu constructed an associative
algebra $A(V)$ [Z] such that there is a one-to-one correspondence 
between irreducible admissible $V$-modules and irreducible $A(V)$-modules.
In the case that $V$ is rational the admissible $V$-module category
and $A(V)$-module category are in fact equivalent. But if $V$ is 
not rational, $A(V)$ does not carry enough information for
the representations of $V.$

In this paper we construct a sequence of associative algebras $A_n(V)$ 
($n=0,1,2,...$) such that $A_0(V)=A(V)$ and $A_n(V)$ is an epimorphic image
of $A_{n+1}(V).$  As in [Z], we use $A_n(V)$ to study representation theory
of $V.$ Let $M=\oplus_{k\geq 0}M(k)$ be an admissible $V$-module
as defined in [DLM] with $M(0)\ne 0.$ Then each $M(k)$ for $k\leq n$ is
an $A_n(V)$-module. In some sense, $A_n(V)$ takes care of the first
$n+1$ homogeneous subspaces of $M$ while $A(V)$ concerns about the top
level $M(0).$ The results of the present paper are modeled on the 
results in [DLM] and the methods are also similar. 
However, the situation for constructing admissible $V$-modules 
from $A_n(V)$-modules turns out to be very complicated. As in [L2] and
[DLM] we extensively use the Lie algebra 
$$\hat V=V\otimes \C[t,t^{-1}]/(L(-1)\otimes 1 +1\otimes {d\over dt})(V\otimes \C[t,t^{-1}])$$
to construct admissible $V$-modules from $A_n(V)$-modules. 
 
It should be pointed out that $\{A_n(V)\}$ in fact form an inverse 
system. So it is natural to consider the inverse limit 
$\displaystyle{\lim_{\leftarrow}A_n(V)}$ and its representations.
This problem will be addressed in a separate paper.  

One of important motivations for constructing $A_n(V)$ is to 
study induced modules from a subalgebra to $V$ as initiated in [DL].
Induced module theory is very important in the representation theory
of classical objects such as groups, rings, Lie algebras. The 
theory of $A_n(V)$ developed in this paper will definitely play 
a role in the study of induced modules for vertex operator algebras.
In order to see this, we consider a subalgebra $U$ of $V$ and 
a $U$-submodule $W$ of $M$ which is an admissible $V$-module. In general,
the top level of $W$ is not necessarily a subspace of the top level of $M.$
In other words, an $A(U)$-module can be a subspace of an $A_n(V)$-module
for some $n>0.$ 
One can now see how $A_n(V)$ enter the picture of studying the 
induced module for the pair $(U,V)$ along this line. 

This paper is organized as follows: In Section 2 we introduce the algebra
$A_n(V)$ which is a quotient of $V$ modulo a subspace $O_n(V)$ consisting
of $u\circ_n v$ (see Section 2 for the definition)
and $L(-1)u+L(0)u$ for $u,v\in V.$ In the case
$n=0,$ $(L(-1)+L(0))u$ can be expressed as $\omega\circ_0 u.$
But in general it is not clear if one can write $(L(-1)+L(0))u$ as 
a linear combination of $v\circ_n w$'s. On the other hand, 
weight zero component of the vertex operator 
$Y((L(-1)+L(0))u,z)$ is zero on any weak $V$-module. 
So we have to put $(L(-1)+L(0))V$ artificially in $O_n(V)$ for general 
$n.$ We also show
in this section how the identity map on $V$ induces an epimorphism
of algebras from $A_{n+1}(V)$ to $A_n(V).$ In Section 3, we construct a
functor $\O_n$ from the category of weak $V$-modules to the category
of $A_n(V)$-modules such that if $M=\oplus_{k\geq 0}M(k)$
is an admissible $V$-module then $\oplus_{k=0}^nM(k)$ with $M(0)\ne 0$ 
is contained  in $\O_n(M)$  and each $M(k)$ for $k\leq n$ is an
$A_n(V)$-submodule. In particular, if $M$ is irreducible then
 $\oplus_{k=0}^nM(k)=\O_n(M)$ and each $M(k)$ is an irreducible
$A_n(V)$-module. 

Section 4 is the core of this paper. In this section we construct a functor
$L_n$ from the category of $A_n(V)$-modules which cannot factor through
$A_{n-1}(V)$ to the category of admissible $V$-modules. For any such
$A_n(V)$-module $U$ we first construct a universal admissible
$V$-module $\bar M_n(U)$ which is somehow a ``generalized Verma module.''
The $L_n(V)$ is then a suitable quotient of $\bar M_n(U);$ the proof 
of this result is technically the most difficult part of this paper.
We also show that $\O_n(L_n(U))/\O_{n-1}(L_n(U))$ is isomorphic to $U$ as 
$A_n(V)$-modules. Moreover, $V$ is rational if and only if
$A_n(V)$ are finite-dimensional semisimple
algebras for all $n.$ Section 5 deals with 
several combinatorial identities used in previous sections.

We assume that the reader is familiar with the basic knowledge on 
vertex operator algebras as presented in  [B],
[FHL] and [FLM]. We also refer the reader to [DLM] for the definitions
of weak modules, admissible modules and (ordinary) modules. 

\section{The associative algebra $A_n(V)$}

Let $V=(V,Y,{\bold 1},\omega)$ be a vertex operator algebra. 
We will construct an associative algebra $A_n(V)$ 
for any nonnegative integer $n$ generalizing the Zhu's algebra $A(V)$
which is our $A_0(V).$

Let $O_n(V)$ be the linear span of all $u\circ_n v$ and $L(-1)u+L(0)u$
where for homogeneous $u\in V$ and $v\in V,$
\bee{g2.2}
u\circ_n v=\Res_{z}Y(u,z)v\frac{(1+z)^{\wt u+n}}{z^{2n+2}}.
\end{equation}
Define the linear space $A_n(V)$ to be the quotient $V/O_{n}(V).$ 

We also define a second product $*_n$ on $V$ for $u$ and $v$ as
above: 
\bee{a5.1}
u*_nv=\sum_{m=0}^{n}(-1)^m{m+n\choose n}\Res_zY(u,z)\frac{(1+z)^{\wt\,u+n}}{z^{n+m+1}}v.
\end{equation} 
Extend linearly to obtain a bilinear product  on $V$ which coincides with that
of Zhu [Z] if $n=0.$ We denote the product (\ref{a5.1}) by $*$
in this case. 
Note that (\ref{a5.1}) may be
written in the form 
\begin{equation}\label{g2.4}
u*_nv=\sum_{m=0}^n\sum_{i=0}^{\infty}(-1)^m
{m+n\choose n}{\wt u+n\choose i}u_{i-m-n-1}v.
\end{equation}

Following lemma generalizes Lemmas 2.1.2 and 2.1.3 of [Z].
\begin{lem}\label{l2.2} (i) Assume that $u\in V$ is homogeneous,
$v\in V$ and $m\ge k\ge 0.$ Then 
$$\Res_{z}Y(u,z)v\frac{(1+z)^{{\wt}u+n+k}}{z^{2n+2+m}}\in O_n(V).$$

(ii) Assume that $v$ is also homogeneous. Then
$$u*_nv-\sum_{m=0}^n{m+n\choose n}(-1)^n\Res_zY(v,z)u\frac{(1+z)^{\wt v+m-1}}{z^{1+m+n}}\in O_n(V)$$
and

(iii) $u*_nv-v*_nu-\Res_zY(u,z)v(1+z)^{\wt u-1}\in O_n(V).$
 \end{lem}

\pf The proof of (i) is similar to that of Lemma 2.1.2 of [Z]. 
As in [Z] we use $L(-1)u+L(0)u\in O_n(V)$ to derive the following formula
$$Y(u,z)v\equiv(1+z)^{-\wt u-\wt v}Y(v,\frac{-z}{1+z})u\ \ \mod O_n(V).$$
Thus we have 
\begin{eqnarray*}
& &u*_nv=\sum_{m=0}^n(-1)^m{n+m\choose n}\Res_zY(u,z)v\frac{(1+z)^{\wt u+n}}{z^{m+n+1}}\\
& &\ \ \ \equiv \sum_{m=0}^n(-1)^m{n+m\choose n}\Res_zY(v,\frac{-z}{1+z})u\frac{(1+z)^{-\wt v+n}}{z^{m+n+1}}\ \ \mod O_n(V)\\
& &\ \ \ =\sum_{m=0}^n(-1)^n{n+m\choose n}\Res_zY(v,z)u\frac{(1+z)^{\wt v+m-1}}{z^{m+n+1}}
\end{eqnarray*}
and (ii) is proved. 

Using (ii) we have 
\begin{eqnarray*}
& &\ \  \ u*_nv-v*_nu\\
& &\equiv \Res_zY(u,z)v(1+z)^{\wt u-1}\left(\sum_{m=0}^n{m+n\choose n}\frac{(-1)^m(1+z)^{n+1}-(-1)^n(1+z)^m}{z^{n+m+1}}\right).
\end{eqnarray*}
By Proposition \ref{ap2} in the Appendix we know that 
$$\sum_{m=0}^n{m+n\choose n}\frac{(-1)^m(1+z)^{n+1}-(-1)^n(1+z)^m}{z^{n+m+1}}=1.$$ 
The proof is complete. \qed

\begin{lem}\label{l2.3} $O_n(V)$ is a 2 sided ideal of $V$ under $*_n.$
\end{lem}

\pf First we show that $(L(-1)u+L(0)u)*_nv\in O_n(V)$ for any homogeneous
$u\in V.$ From the definition we see that
\begin{eqnarray*}
& &(L(-1)u)*_nv=\sum_{m=0}^{n}{m+n\choose n}(-1)^m\Res_zY(L(-1)u,z)
\frac{(1+z)^{\wt\,u+n+1}}{z^{n+m+1}}v\\
& &\ \ \ =\sum_{m=0}^{n}{m+n\choose n}(-1)^m\Res_z({d \over dz}Y(u,z))v
\frac{(1+z)^{\wt\,u+n+1}}{z^{n+m+1}}\\
& &\ \ \ =\sum_{m=0}^{n}{m+n\choose n}(-1)^{m+1}\Res_zY(u,z)v\left(
\frac{(-n-m-1)(1+z)^{\wt\,u+n+1}}{z^{n+m+2}}\right.\\
& &\ \ \ \ \ \ \ \ \left.+\frac{z(\wt u+n+1)(1+z)^{\wt\,u+n}}{z^{n+m+2}}\right).
\end{eqnarray*}
Thus 
$$(L(-1)u+\wt u\,u)*_nv=\sum_{m=0}^{n}{m+n\choose n}(-1)^{m}
\Res_zY(u,z)v(1+z)^{\wt u+n}\frac{mz+n+m+1}{z^{n+m+2}}.$$
It is straightforward to show that 
\begin{eqnarray*}
& &\ \ \  \sum_{m=0}^{n}{m+n\choose n}(-1)^{m}\frac{mz+n+m+1}{z^{n+m+2}}\\
& &=\sum_{m=0}^{n}{m+n\choose n}(-1)^{m}\frac{mz}{z^{n+m+2}}
+\sum_{m=0}^{n}{m+n+1\choose n}(-1)^{m}\frac{m+1}{z^{n+m+2}}\\
& &=(-1)^n{2n\choose n}\frac{2n+1}{z^{2n+2}}.
\end{eqnarray*}
It is clear now that $(L(-1)u+L(0)u)*_nv\in O_n(V).$

Second, we show that $u*_n(L(-1)v+L(0)v)\in O_n(V).$ Using the result
that  $(L(-1)v+L(0)v)*_nv\in O_n(V)$ and Lemma \ref{l2.2} (iii)
we have
\begin{eqnarray*}
& &\ \ u*_n(L(-1)v+L(0)v)\\
& &\equiv -\Res_z\left(Y(L(-1)v,z)u(1+z)^{\wt v}+Y(L(0)v,z)u(1+z)^{\wt v-1}\right)\ \ \mod O_n(V)\\
& &=\Res_z\left(Y(v,z)u\frac{d}{dz}(1+z)^{\wt v}-Y(L(0)v,z)u(1+z)^{\wt v-1}\right)\\
& &=0.
\end{eqnarray*}

Third, a similar argument as in [Z] using Lemma \ref{l2.2} (i) shows that
$u*_n(v\circ_n w)\in O_n(V)$ for $u,v,w\in V.$

Finally, use $u*_n(v\circ_n w)\in O_n(V)$ and 
Lemma \ref{l2.2} (iii) to obtain
\begin{eqnarray*}
& &\ \ \ (v\circ_n w)*_n u\\
& &\equiv -\Res_{z_1}\Res_{z_2}Y(u,z_1)Y(v,z_2)w
\frac{(1+z_1)^{\wt u-1}(1+z_2)^{\wt v+n}}{z_2^{2n+2}}\ \ \mod O_n(V)\\
& &\equiv -\Res_{z_2}\Res_{z_1-z_2}Y(Y(u,z_1-z_2)v,z_2)w
\frac{(1+z_1)^{\wt u-1}(1+z_2)^{\wt v+n}}{z_2^{2n+2}}\\
& &=-\sum_{i\geq 0}{\wt u-1\choose i}\Res_{z_2}Y(u_iv,z_2)w\frac{(1+z_2)^{\wt u+\wt v+n-1-i}}{z_2^{2n+2}}
\end{eqnarray*}
which belongs to $O_n(V)$ as $\wt u_iv=\wt u+\wt v-i-1.$ This completes
the proof.
\qed

Our first main result is the following.
\begin{thm}\label{t2.4}  (i) The product $*_n$ induces the structure of an 
associative algebra  on $A_n(V)$ with identity ${\bf 1}+O_n(V).$

(ii) The linear map 
$$\phi:  v\mapsto e^{L(1)}(-1)^{L(0)}v$$
induces an anti-isomorphism $A_n(V)\to A_n(V)$.

(iii) $\omega+O_n(V)$ is a central element of $A_n(V).$
\end{thm}

\pf For (i) we only need to prove that $A_n(V)$ is associative.
Let $u,v,w\in V$ be homogeneous. Then
\begin{eqnarray*}
& &(u*_n v)*_nw=\sum_{m_1=0}^n\sum_{i\geq 0}(-1)^{m_1}{m_1+n\choose
n}{\wt u+n\choose i}(u_{-m_1-n-1+i}v)*_nw\\
& &=\sum_{m_1,m_2=0}^n\sum_{i\geq 0}(-1)^{m_1+m_2}{m_1+n\choose
n}{m_2+n\choose n}{\wt u+n\choose i}\\
& &\ \ \ \ \ \Res_{z_2}Y(u_{-m_1-n-1+i}v,z_2)w
\frac{(1+z_2)^{\wt u+\wt v+2n+m_1-i}}{z_2^{1+m_2+n}}\\
& &=\sum_{m_1,m_2=0}^n(-1)^{m_1+m_2}{m_1+n\choose
n}{m_2+n\choose n}\\
& &\ \ \ \ \ \Res_{z_2}\Res_{z_1-z_2}Y(Y(u,z_1-z_2)v,z_2)w 
\frac{(1+z_1)^{\wt u+n}(1+z_2)^{\wt v+n+m_1}}{(z_1-z_2)^{m_1+n+1}z_2^{1+m_2+n}}\\
& &=\sum_{m_1,m_2=0}^n(-1)^{m_1+m_2}{m_1+n\choose
n}{m_2+n\choose n}\\
& &\ \ \ \ \ \Res_{z_1}\Res_{z_2}Y(u,z_1)Y(v,z_2)w
\frac{(1+z_1)^{\wt u+n}(1+z_2)^{\wt v+n+m_1}}{(z_1-z_2)^{m_1+n+1}z_2^{1+m_2+n}}\\
& &\ \ \ \ \ -\sum_{m_1,m_2=0}^n(-1)^{m_1+m_2}{m_1+n\choose
n}{m_2+n\choose n}\\
& &\ \ \ \ \ \Res_{z_2}\Res_{z_1}Y(v,z_2)Y(u,z_1)w\frac{(1+z_1)^{\wt u+n}(1+z_2)^{\wt v+n+m_1}}{(z_1-z_2)^{m_1+n+1}z_2^{1+m_2+n}}\\
& &=\sum_{m_1,m_2=0}^n\sum_{i\geq 0}(-1)^{m_1+m_2}{m_1+n\choose
n}{m_2+n\choose n}{-m_1-n-1\choose i}(-1)^i\\
& &\ \ \ \ \ \Res_{z_1}\Res_{z_2}Y(u,z_1)Y(v,z_2)w\frac{(1+z_1)^{\wt u+n}(1+z_2)^{\wt v+n+m_1}}{z_1^{m_1+n+1+i}z_2^{1+m_2+n-i}}\\
& &\ \ \ \ \ -\sum_{m_1,m_2=0}^n\sum_{i\geq 0}(-1)^{m_2+n+1+i}{m_1+n\choose
n}{m_2+n\choose n}{-m_1-n-1\choose i}\\
& &\ \ \ \ \ \Res_{z_2}\Res_{z_1}Y(v,z_2)Y(u,z_1)w
\frac{(1+z_1)^{\wt u+n}(1+z_2)^{\wt v+n+m_1}}{z_1^{-i}
z_2^{2+m_1+m_2+2n+i}}.
\end{eqnarray*}
From Lemma \ref{l2.2} we know that 
$$\Res_{z_2}\Res_{z_1}Y(v,z_2)Y(u,z_1)w
\frac{(1+z_1)^{\wt u+n}(1+z_2)^{\wt v+n+m_1}}{z_1^{m_1+n+1-i}z_2^{2+m_1+m_2+2n+i}}$$
lies in $O_n(V).$ Also if $i>n-m_1$ 
$$\Res_{z_1}\Res_{z_2}Y(u,z_1)Y(v,z_2)w\frac{(1+z_1)^{\wt u+n}(1+z_2)^{\wt v+n+m_1}}{z_1^{m_1+n+1+i}z_2^{1+m_2+n-i}}$$
is in $O_n(V).$ Thus 
\begin{eqnarray*}
& &(u*_n v)*_nw\equiv u*_n(v*_n w)+\sum_{m_1,m_2=0}^n(-1)^{m_1+m_2}{m_1+n\choose
n}{m_2+n\choose n}\\
& &\ \ \ \ \ \Res_{z_1}\Res_{z_2}Y(u,z_1)Y(v,z_2)\frac{(1+z_1)^{\wt u+n}(1+z_2)^{\wt v+n}}{z_1^{m_1+n+1}z_2^{1+m_2+n}}\\
& &\ \ \ \ \ \left(\sum_{i=0}^{n-m_1}\sum_{j\geq 0}
{-m_1-n-1\choose i}{m_1\choose j}(-1)^i\frac{z_2^{i+j}}{z_1^i}-1\right).
\end{eqnarray*}
From Proposition \ref{ap3} in the Appendix we know that 
$$\sum_{m_1=0}^n(-1)^{m_1}{m_1+n\choose
n}\left(\sum_{i=0}^{n-m_1}\sum_{j\geq 0}
{-m_1-n-1\choose i}{m_1\choose j}(-1)^i\frac{z_2^{i+j}}{z_1^{i+m_1}}-1\right)=0.$$
This implies that the product $*_n$ of $A_n(V)$ is associative.

The proof of (ii) is similar to that of (ii) of Theorem 2.4 [DLM]. 
We refer the reader to [DLM] for detail.

Note that ${\bf 1}*_n u=u$ for any $u\in V$ and that
$$u*_n {\bf 1}-{\bf 1}*_n u\equiv \Res_zY(u,z){\bf 1}(1+z)^{\wt u-1}=0.$$
This shows that ${\bf 1}+O_n(V)$ is the identity of $A_n(V).$
Again by Lemma \ref{l2.2} (iii),
$$\omega*_n u-u*_n\omega=\Res_zY(\omega,z)u(1+z)=L(-1)u+L(0)u\in O_n(V).$$
So (iii) is proved. \qed

\begin{prop}\label{inver}
The identity map on $V$ induces an onto algebra homomorphism
from $A_n(V)$ to $A_{n-1}(V).$
\end{prop}

\pf First by Lemma \ref{l2.2} (i), $O_n(V)\subset O_{n-1}(V).$ It remains to
show that $u*_nv\equiv u*_{n-1}v \ \mod O_{n-1}(V).$ Let $u$ be homogeneous.
Then 
\begin{eqnarray*}
& &u*_nv=\sum_{m=0}^{n}{m+n\choose n}(-1)^m \Res_zY(u,z)v\frac{(1+z)^{{\wt}\,u+n-1}}{z^{n+m}}\\
& & \ \ \ \ \ +\sum_{m=0}^{n}{m+n\choose n}(-1)^m \Res_zY(u,z)v\frac{(1+z)^{{\wt}\,u+n-1}}{z^{n+m+1}}\\
& &\equiv \sum_{m=0}^{n-1}{m+n\choose n}(-1)^m\Res_zY(u,z)v\frac{(1+z)^{{\wt}\,u+n-1}}{z^{n+m}}\ \ \\
& &\ \ \ \ \  +\sum_{m=0}^{n-2}{m+n\choose n}(-1)^m\Res_zY(u,z)\frac{(1+z)^{{\wt}\,u+n-1}}{z^{n+m+1}}\ \mod O_{n-1}(V)\\
& &=\Res_zY(u,z)v\frac{(1+z)^{\wt u+n-1}}{z^n}+\!\sum_{m=1}^{n-1}\Res_zY(u,z)v\frac{(1\!+\!z)^{{\wt}\,u+n-1}}{z^{n+m}}\cdot\\
& &\ \ \ \ \ \cdot\left((-1)^m{m+n\choose n}\!+\!(-1)^{m+1}{m+n-1\choose n}\right)\\
& &=u*_{n-1}v,
\end{eqnarray*}
as desired.
\qed

From Proposition \ref{inver} we in fact have an inverse system $\{A_n(V)\}.$
Denote by $I(V)$ the inverse limit $\displaystyle{\lim_{\leftarrow}A_n(V).}$
Then\begin{equation}\label{i1}
I(V)=\{a=(a_n+O_n(V))\in \prod_{n=0}^{\infty}A_n(V)| 
a_n-a_{n-1}\in O_{n-1}(v)\}.
\end{equation}
Define $i:V\to I(V)$ such that $i(v)=(v+O_n(V))$ for $v\in V.$ Then $V/\ker i$
is linearly isomorphic to a subspace of $I(V).$ It is easy to see that $i(V)$
is not closed under the product. But one can introduce a suitable
topology on $I(V)$ so that $i(V)$ is a dense subspace of $I(V)$ under
the topology. An interesting problem is to determine the kernel
of $i. $ From the definition of $O_n(V)$ we see immediately that
$(L(-1)+L(0))V$ is contained in the kernel. It will be proved in Section
3 that if $v\in O_n(V)$ then
 $a_{\wt v-1}=0$ on $\oplus_{k=0}^nM(n)$ for any admissible $V$-module
$\oplus_{k=0}^{\infty}M(k).$  Thus $a\in \ker i$ if and only if
$a_{\wt a-1}=0$ on any admissible $V$-module. It is proved in
[DLMM] that if $V$ is a simple vertex operator algebra $V$ satisfies 
$V_k=0$ for
$k<0$ and $V_0=\C{\bf 1}$ then the subspace of $V$ consisting
of vectors $v$ whose component operators $v_{\wt v-1}$ are 0  on $V$ 
is essentially $(L(0)+L(-1))V.$ We suspect that if $V$ is a rational
vertex operator algebra then the kernel of $i$ is exactly $(L(0)+L(-1))V.$

\section{The Functor $\Omega_n$}

Consider the quotient space
\begin{equation}\label{g4.1}
\hat V={\C}[t,t^{-}]\otimes V/D{\C}[t,t^{-}]\otimes V
\end{equation}
where $D={d\over dt}\otimes 1+1\otimes L(-1).$ Denote by $v(m)$ the image
of $v\otimes t^m$ in $\hat V$ for $v\in V$ and $m\in \Z.$ Then $\hat V$
is $\Z$-graded by defining the degree of $v(m)$ to be $\wt v-m-1$ if $v$
is homogeneous. Denote the homogeneous subspace of degree $m$ by $\hat V(m).$
The space $\hat V$ is, in fact, a $\Z$-graded Lie algebra with bracket
\begin{equation}\label{ld}
[a(p), b(q)]=\sum_{i=0}^{\infty}{p\choose i}a_ib(p+q-i)
\end{equation}
(see [L2] and [DLM]). In particular, $\hat V(0)$ is a Lie subalgebra. 
By Lemma \ref{l2.2} (iii) we have 
\bp{p2.5}
 Regarded $A_n(V)$ as a Lie algebra, the map $v(\wt v-1)\mapsto v+O_n(V)$ 
is a well-defined onto Lie algebra homomorphism from $\hat V(0)$ to
$A_n(V).$
\ep

By Lemmas 5.1 and 5.2 of [DLM], any weak $V$-module
 $M$ is a module for $\hat V$ under the map map $a(m)\mapsto a_m$ and
a weak $V$-module which carries a 
$\Z_+$-grading is an admissible 
$V$-module if, and only if, $M$ is a $\Z_+$-graded module for
the graded Lie algebra $\hat V.$

For a module $W$ for the Lie algebra $\hat V$ and a nonnegative
$m$ we let $\O_m(W)$ denote the space
of ``m-th lowest
weight vectors,'' that is
\begin{eqnarray}\label{g5.1}
\Omega_m(W)=\{u\in W|\hat V(-k)u=0\ if \ k\geq m\}.
\end{eqnarray}
Then $\O_m(W)$ is a module for the Lie algebra $\hat V(0).$

\bt{t5.3} Suppose that $M$ is a weak $V$-module. Then there is 
a representation of the associative algebra $A_n(V)$ on $\O_n(M)$ induced by
the map $a\mapsto o(a)=a_{\wt a-1}$ for homogeneous $a\in V.$ 
\et

\pf We need to show that $o(a)=0$ for all $a\in O_n(V)$ and
$o(u*_nv)=o(u)o(v)$ for $u,v\in V.$
Using $Y(L(-1)u,z)=\frac{d}{dz}Y(u,z)$ we immediately see that 
$o(L(-1)u+L(0)u)=0.$ From the proof of Lemma \ref{l2.2} we know that
$(L(-1)u+L(0)u)*_nv=(-1)^n{2n\choose n}(2n+1)u\circ_n v.$ It suffices
to show that $o(u*_nv)=o(u)o(v).$ 

Let $u,v$ be homogeneous and $0\leq k\leq n.$ 
Note that $v_{\wt v+p}=u_{\wt u+p}=0$ on $\Omega_n(M)$ 
if $p\geq n.$ We assert that the following identity holds on $\Omega_n(M):$
\begin{eqnarray}\label{au1}
& & \sum_{m=0}^k(-1)^m{2n+m-k\choose m}o(\Res_zY(u,z)v\frac{(1+z)^{\wt u+n}}{z^{2n+1-k+m}})\nonumber\\
& &\ \ \ \ =u_{\wt u-n+k-1}v_{\wt v+n-k-1} 
\end{eqnarray}
which reduces to $o(u*_n v)=o(u)o(v)$ if $k=n.$  
The proof of (\ref{au1}) is a straightforward computation involving the 
Jacobi identity on modules in terms of residues. 

On $\Omega_n(M)$ we have
\begin{eqnarray*}
& &\sum_{m=0}^k(-1)^m{2n+m-k\choose m}o(\Res_zY(u,z)v\frac{(1+z)^{\wt u+n}}{z^{2n+1-k+m}})\\
& &=\sum_{m=0}^k\sum_{i\geq 0}(-1)^m{2n+m-k\choose m}{\wt u+n\choose i}
o(u_{i-2n-1-m+k}v)\\
& &=\sum_{m=0}^k\sum_{i\geq 0}(-1)^m{2n+m-k\choose m}{\wt u+n\choose i}
 (u_{i-2n-1-m+k}v)_{\wt u+\wt v-i+2n+m-1-k}\\
& &=\sum_{m=0}^n\sum_{i\geq 0}(-1)^m{2n+m-k\choose m}{\wt u+n\choose i}\\
& &\ \ \ \ \Res_{z_2}\Res_{z_1-z_2}Y(Y(u,z_1-z_2)v,z_2)(z_1-z_2)^{i-2n-m-1+k}z_2^{\wt u+\wt v-i+2n+m-1-k}\\
& &=\sum_{m=0}^k(-1)^m{2n+m-k\choose m}\Res_{z_2}\Res_{z_1-z_2}Y(Y(u,z_1-z_2)v,z_2)\frac{z_1^{\wt u+n}z_2^{\wt v+n+m-1-k}}{(z_1-z_2)^{2n+m+1-k}}\\
& &=\sum_{m=0}^k(-1)^m{2n+m-k\choose m}\Res_{z_1}\Res_{z_2}Y(u,z_1)Y(v,z_2)\frac{z_1^{\wt u+n}z_2^{\wt v+n+m-1-k}}{(z_1-z_2)^{2n+m+1-k}}\\
& &-\sum_{m=0}^k(-1)^m{2n+m-k\choose m}\Res_{z_2}\Res_{z_1}Y(v,z_2)Y(u,z_1)\frac{z_1^{\wt u+n}z_2^{\wt v+n+m-1-k}}{(z_1-z_2)^{2n+m+1-k}}\\
& &=\sum_{m=0}^k\sum_{i=0}^{k\!-\!m}(-1)^{m+i}{2n\!+\!m\!-\!k\choose m}{-\!m\!-\!2n\!-1\!+\!k\choose i}u_{\wt u-n-m-1+k-i}v_{\wt v+n+m-1-k+i}\\
& &=\sum_{m=0}^{k}\sum_{i=m}^{k}{2n+m-k\choose m}{-m-2n-1+k\choose i-m}(-1)^{i}u_{\wt u-n+k-i-1}v_{\wt v+n-k-1+i}\\
& &=\sum_{i=0}^{k}\sum_{m=0}^{i}{2n+m-k\choose m}{-m-2n-1+k\choose i-m}(-1)^{i}u_{\wt u-n+k-i-1}v_{\wt v+n-k-1+i}\\
& &=u_{\wt u-n+k-1}v_{\wt v+n-k-1}\\
& &+\sum_{i=1}^{k}\sum_{m=0}^{i}{2n+m-k\choose m}{-m-2n-1+k\choose i-m}(-1)^{i}u_{\wt u-n+k-i-1}v_{\wt v+n-k-1+i}.
\end{eqnarray*}
It is enough to show that for $i=1,...,k,$
$$\sum_{m=0}^{i}{2n+m-k\choose m}{-m-2n-1+k\choose i-m}=0,$$
which follows from a easy calculation:
\begin{eqnarray*}
& &\ \ \sum_{m=0}^{i}{2n+m-k\choose m}{-m-2n-1+k\choose i-m}\\
& &=\sum_{m=0}^{i}(-1)^{i-m}{2n+m-k\choose m}{2n+i-k\choose i-m}\\
& &=\sum_{m=0}^{i}(-1)^{i-m}{2n+i-k\choose 2n-k}{i\choose m}\\
& &=0.
\end{eqnarray*}
This completes the proof. \qed

\begin{rem}\label{rim} For homogeneous $u,v\in V$ and $j\in\Z$ we set 
$o_j(u)=u_{\wt u-1-j}$ and extend to all $u\in V$ by linearity. 
Then $o_0(u)=o(u).$ Using associativity of the vertex operators
 $$(z_0+z_2)^{\wt u+n}Y(u,z_0+z_2)Y(v,z_2)=(z_2+z_0)^{\wt u+n}Y(Y(u,z_0)v,z_2)$$
on $\Omega_n(M)$ we have that for $i\geq j$ with $i+j\geq 0$ 
these exists a unique $w_{u,v}^{i,j}\in V$ such that 
$o_i(u)o_j(v)=o_{i+j}(w_{u,v}^{i,j})$ on  $\Omega_n(M).$ 
In fact one can write down $w_{u,v}^{i,j}$ explicitly in terms of 
$u$ and $v.$ But for our later purpose it is enough to know
the explicit expression of $w_{u,v}^{i,-i}$ ($i\geq 0$) which
is given by
$$w_{u,v}^{i,-i}=\sum_{m=0}^{n-i}(-1)^m{n+m+i\choose m}\Res_zY(u,z)v\frac{(1+z)^{\wt u+n}}{z^{n+1+i+m}}$$
in the proof of Theorem \ref{t5.3}.
\end{rem}  

It is clear that $\Omega_n$ is
a covariant functor from the category of 
weak $V$-modules 
to the category of $A_n(V)$-modules. To be more precise, if $f:M\to N$
is a morphism in the first category  we define $\O_n(f)$
to be the restriction of $f$ to $\O_n(M)$. 
Then $f$ induces a morphism of $\hat V$-modules $M\to N$
by Lemma 5.1 of [DLM]. Moreover
$\O_n(f)$ maps $\O_n(M)$ to $\O_n(N).$ Now Theorem \ref{t5.3} implies that $\O_n(f)$ is a morphism of $A_n(V)$-modules.

Let $M$ be such a module. As long as $M\ne 0,$ then some
$M(m)\ne 0,$  and it is no loss to shift the grading
so that in fact $M(0)\ne 0.$ If $M=0,$ let $M(0)=0.$ With these
conventions we prove 
\bp{l2.9} Suppose that $M$ is an admissible $V$-module. 
Then the following hold

(i) $\Omega_n(M)\supset \oplus_{i=0}^{n}M(i).$ If $M$ is simple
then  $\Omega_n(M)=\oplus_{i=0}^{n}M(i).$

(ii) Each $M(p)$ is an $\hat V(0)$-module and
$M(p)$ and $M(q)$ are inequivalent if $p\ne q$ and both
$M(p)$ and $M(q)$ are nonzero.  If $M$ is simple then each
$M(p)$ is an irreducible $\hat V(0)$-module.

(iii) Assume that $M$ is simple. Then
each $M(i)$ for $i=0,...,n$ is a simple $A_n(V)$-module
and $M(i)$ and $M(j)$ are inequivalent $A_n(V)$-modules. 
\ep

\pf An easy argument shows that $\O_n(M)$ is a graded subspace of $M.$ That is
\begin{equation}\label{g5.6}
\O_n(M)=\oplus_{i\geq 0}\O_n(M)\cap M(i).
\end{equation}
Set $\O_n(i)=\O(M)\cap M(i).$ It is clear that $M(i)\subset \O_n(M)$ if
$i<n.$ In order to prove (i) we must show that $\O_n(i)=0$ if $i\geq n.$

By Proposition 2.4 of [DM] or Lemma 6.1.1 of [L2], $M=\spa\{u_nw|u\in
V,n\in\Z\}$ where $w$ is any nonzero vector in $M.$ If $\O_n(i)\ne 0$ for
some $i\geq n$ we can take $0\ne w\in
\O_n(i).$ Since $u_{\wt u+p}w=0$ for all $p\geq n$
we see that $M=\spa\{u_{\wt u+p}w|u\in V,p\in\Z, p<n\}.$ This implies
that $M(0)=0, $ a contradiction. 

It is clear that (iii) follows from (ii). For (ii), note that
$M=\hat Vw=\oplus_{p\in\Z}\hat V(p)w.$ Thus if $0\ne w\in M(i)$ 
then $\hat V(p)w=M(i+p).$ In particular, $\hat V(0)w=M(i),$
as required. It was pointed out in [Z] that $L(0)$ is semisimple
on $M$ and $M(k)=\{w\in M| L(0)w=(h+k)w\}$ for some fixed $h.$ 
The inequivalence follows from the fact that $L(0)$
has different eigenvalues on $M(p)$ and $M(q).$
\qed

\section{The functor $L_n$}

We show in this section that there is a universal
way to construct an admissible $V$-module
from an $A_n(V)$-module which cannot factor through $A_{n-1}.$
(If it can factor through $A_{n-1}(V)$ we can consider the same
procedure for $A_{n-1}(V).$) 
Moreover a certain quotient of the universal object is an
admissible $V$-module $L_n(U)$ and $L_n$ defines a functor which
is a right inverse to the functor $\O_n/\O_{n-1}$ 
where $\O_n/\O_{n-1}$ is the quotient functor $M\mapsto \O_n(M)/\O_{n-1}(M).$

Fix an $A_n(V)$-module $U$ which cannot factor through $A_{n-1}(V).$ Then
it is a  module for $A_n(V)_{Lie}$ in an obvious way.
By Proposition \ref{p2.5} we can lift $U$ to a module for the Lie algebra
$\hat V(0),$ and then to one for $P_n=\oplus_{p>n}\hat V(-p)\oplus \hat V(0)$ by letting $\hat V(-p)$ act trivially. Define 
\bee{g6.1}
M_n(U)=\Ind_{P_n}^{\hat V}(U)=U(\hat V)\otimes_{U(P_n)} U. 
\end{equation}
If we give $U$ degree $n$, the $\Z$-gradation of $\hat V$ lifts to 
$M_n(U)$ which thus becomes a  $\Z$-graded module for $\hat V.$
The PBW theorem implies that $M_n(U)(i)=U(\hat V)_{i-n}U.$ 

We define for $v\in V,$
\bee{g6.2}
Y_{M_n(U)}(v,z)=\sum_{n\in\Z}v(m)z^{-m-1}
\end{equation}
As in [DLM], $Y_{M(U)}(v,z)$ satisfies all conditions of a week $V$-module
except the associativity which does not hold on $M_n(U)$ in general.
We have to divide out by the desired relations. 

Let $W$ be the subspace of $M_n(U)$ spanned linearly by the 
coefficients of
\begin{eqnarray}\label{g6.3}
(z_{0}+z_{2})^{{\wt}a+n}Y(a,z_{0}+z_{2})Y(b,z_{2})u-(z_{2}+z_{0})^{{\wt}a+n}
Y(Y(a,z_{0})b,z_{2})u
\end{eqnarray}
for any homogeneous $a\in V,b\in V,$ $u\in U$.
Set
\be{g6.4}
\bar M_n(U)=M_n(U)/U(\hat V)W.
\end{equation}

\bt{t6.1} The space $\bar M_n(U)$ is an admissible $V$-module $\bar M_n(U)
=\sum_{m\geq 0}\bar M_n(U)(m)$ with $\bar M_n(U)(0)\ne 0,$
$\bar M_n(U)(n)=U$ and with the 
following universal property: for any weak $V$-module $M$
and any $A_n(V)$-morphism $\phi: U\to \O_n(M),$ there is a unique morphism 
$\bar\phi: \bar M_n(U)\to M$ of weak $V$-modules which 
extends $\phi.$ 
\et

\pf By Proposition 6.1 of [DLM], we know that $\bar M_n(U)$ is a $\Z$-graded
weak $V$-module generated by $U+U(\hat V)W.$ By Proposition 2.4 of [DM]
or Lemma 6.1.1 of [L2] $\bar M_n(U)$ is spanned by
$$\{a_n(U+U(\hat V)W)|a\in V,n\in \Z\}.$$
Thus $\bar M_n(U)(m)=\hat V(m-n)(U+U(\hat V)W)$ for all $m\in \Z.$ In 
particular, $\bar M_n(U)(m)=0$ if $m<0$ and 
$\bar M_n(U)(n)=A_n(V)(U+U(\hat V)W)$ which is a quotient module of $U.$
A proof that $\bar M_n(U)(0)\ne 0$ and $\bar M_n(U)(n)=U$
will be given after Proposition \ref{p3.5}.
The universal property of 
$\bar M_n(U)$ follows from its construction.  \qed

In the following we let $U^*=\Hom_{\C}(U,\C)$ and let $U_s$ be the subspace
of $M_n(U)(n)$ spanned by  ``length'' $s$ vectors  
 $$o_{p_1}(a_1)\cdots o_{p_s}(a_s)U$$
where $p_1\geq \cdots \geq p_s,$ $p_1+\cdots p_s=0,$ $p_i\ne 0,$ $p_s\geq -n$
and $a_i\in V.$ Then by PBW theorem 
$M_n(U)(n)=\sum_{s\geq 0}U_s$ with $U_0=U$ and $U_s\cap U_t=0$ if $s\ne t.$
Recall Remark \ref{rim}. We extend $U^*$ to 
$M_n(U)(n)$ inductively so that 
\begin{eqnarray}\label{def}
\<u',o_{p_1}(a_1)\cdots o_{p_s}(a_s)u\>
=\<u',o_{p_{1}+p_2}(w_{a_1,a_2}^{p_1,p_2})o_{p_3}(a_3)\cdots o_{p_{s}}(a_{s})u).
\end{eqnarray}
where $o_j(a)=a(\wt a-1-j)$ for homogeneous $a\in V.$ 
We further extend
$U^*$ to $M_n(U)$ by letting $U^*$ annihilate $\oplus_{i\ne n}M(U)(i).$

Set 
$$ J=\{v\in M_n(U)|\langle u',xv\rangle=0\ {\rm for\ all}\ u'\in
U^{*},\ {\rm all}\ x\in U(\hat V)\}.$$
We can now state the second main result of this section.
\bt{t6.3} Space $L_n(U)\!=\!M_n(U)/J$ is an admissible
$V$-module satisfying $L_n(U)(0)\ne 0$ and
$\O_n/\O_{n-1}(L_n(U))\cong U.$
Moreover 
$L_n$ defines a functor from the category of  
$A_n(V)$-modules which cannot factor through $A_{n-1}(V)$  to
the category of admissible $V$-modules such that $\O_n/\O_{n-1}\circ L_n$
is naturally equivalent to the identity.
\et

The main point in the proof of the Theorem is to show that 
$U(\hat V)W\subset J.$ The next three results are devoted to this goal.

\begin{prop}\label{p3.3} The following hold 
for all homogeneous $a\in V, b\in V,$ $u'\in U^{*}, u\in U,j\in {\Z}_{+},$
\begin{eqnarray}\label{ea3}
& &\ \ \ \langle u',(z_{0}+z_{2})^{{\wt}a+n+j}Y_{M_n(U)}(a,z_{0}+z_{2})Y_{M_n(U)}(b,z_{2})u\rangle
\nonumber\\
& &=\langle u',(z_{2}+z_{0})^{{\wt}a+n+j}Y_{M_n(U)}(Y(a,z_{0})b,z_{2})u\rangle.
\end{eqnarray}
\ep
In the following we simply write $Y$ for $Y_{M_n(U)},$ which
should cause no confusion. The following is the key lemma.
\bl{l6.7} For any $i,j\in \Z_{+},$ 
\begin{eqnarray*}
& &\ \ \Res_{z_{0}}z_{0}^{-1+i}(z_{0}+z_{2})^{{\wt}a+n+j}
\langle u',Y(a,z_{0}+z_{2})Y(b,z_{2})u\rangle\nonumber\\
& &=\Res_{z_{0}}z_{0}^{-1+i}(z_{2}+z_{0})^{{\wt}a+n+j}
\langle u',Y(Y(a,z_{0})b,z_{2})u\rangle.
\end{eqnarray*}
\el

\pf Since $j\geq 0$ then $a(\wt a+n+j)$ lies in $\oplus_{p>n}\hat V(-p)$ and hence
annihilates $u.$ Then for all $i\in\Z_+$ we get
\begin{eqnarray}\label{3.5}
\Res_{z_{1}}(z_{1}-z_{2})^{i}z_{1}^{{\wt}a+n+j}
Y(b,z_{2})Y(a,z_{1})u=0.
\end{eqnarray}
Note that (\ref{ld}) is equivalent to
\begin{eqnarray}\label{3.6''}
[Y(a,z_{1}),Y(b,z_{2})]
=\Res_{z_{0}}z_2^{-1}
\delta\left(\frac{z_1-z_0}{z_2}\right)Y(Y(u,z_0)v,z_2).
\end{eqnarray}

Using (\ref{3.5}) and (\ref{3.6''}) 
we obtain:
\begin{eqnarray}
& &\ \ \ \Res_{z_{0}}z_{0}^{i}(z_{0}+z_{2})^{{\wt}a+n+j}
Y(a,z_{0}+z_{2})Y(b,z_{2})u\nonumber\\
& &= \Res_{z_{1}}(z_{1}-z_{2})^{i}z_{1}^{{\wt}a+n+j}
Y(a,z_{1})Y(b,z_{2})u\nonumber\\
& &= \Res_{z_{1}}(z_{1}-z_{2})^{i}z_{1}^{{\wt}a+n+j}
Y(a,z_{1})Y(b,z_{2})u\nonumber\\
& &\ \ \ -\Res_{z_{1}}(z_{1}-z_{2})^{i}z_{1}^{{\wt}a+n+j}
Y(b,z_{2})Y(a,z_{1})u\nonumber\\
& &= \Res_{z_{1}}(z_{1}-z_{2})^{i}z_{1}^{{\wt}a+n+j}
[Y(a,z_{1}),Y(b,z_{2})]u\nonumber\\
& &= \Res_{z_{0}}\Res_{z_{1}}(z_{1}-z_{2})^{i}z_{1}^{{\wt}a+n+j}
z_{2}^{-1}\delta\left(\frac{z_{1}-z_{0}}{z_{2}}\right)
Y(Y(a,z_{0})b,z_{2})u\nonumber\\
& &=\Res_{z_{0}}\Res_{z_{1}}z_{0}^{i}z_{1}^{{\wt}a+n+j}
z_{1}^{-1}\delta\left(\frac{z_{2}+z_{0}}{z_{1}}\right)
Y(Y(a,z_{0})b,z_{2})u
\nonumber\\
& &=\Res_{z_{0}}z_{0}^{i}(z_{2}+z_{0})^{{\wt}a+n+j}
Y(Y(a,z_{0})b,z_{2})u.\label{3.6}
\end{eqnarray}
Thus lemma \ref{l6.7} holds if $i\geq 1,$ and we may now assume
$i=0.$

Next use (\ref{3.6}) to calculate that 
\begin{eqnarray}
& &\Res_{z_{0}}z_{0}^{-1}(z_{0}+z_{2})^{{\wt}a+n+j}
\langle u',Y(a,z_{0}+z_{2})Y(b,z_{2})u\rangle\nonumber\\
&=&\sum_{k=0}^{\infty}\left(\begin{array}{c}j\\k\end{array}\right)
\Res_{z_{0}}z_{0}^{k-1}z_{2}^{j-k}(z_{0}+z_{2})^{{\wt}a+n}
\langle u',Y(a,z_{0}+z_{2})Y(b,z_{2})u\rangle\nonumber\\
&=&\sum_{k=1}^{\infty}\left(\begin{array}{c}j\\k\end{array}\right)
\Res_{z_{0}}z_{0}^{k-1}z_{2}^{j-k}(z_{2}+z_{0})^{{\wt}a+n}
\langle u',Y(Y(a,z_{0})b,z_{2})u\rangle\nonumber\\
& & \ \ \ + \Res_{z_{0}}z_{0}^{-1}z_2^j(z_{2}+z_{0})^{{\wt}a+n}
\langle u',Y(a,z_0+z_2)Y(b,z_{2})u\rangle.\label{g6.9}
\end{eqnarray}

It reduces to show that 
\begin{eqnarray}\label{g6.10}
& & \Res_{z_{0}}z_{0}^{-1}(z_{2}+z_{0})^{{\wt}a+n}
\langle u',Y(a,z_0+z_2)Y(b,z_{2})u\rangle\\
& &= \Res_{z_{0}}z_{0}^{-1}(z_{2}+z_{0})^{{\wt}a+n}
\langle u',Y(a,z_0)u,z_2)u\rangle.
\end{eqnarray}

Since $\<u',M_n(U)(m)\>=0$ if $m\ne n,$  we see that
\begin{eqnarray*}
& &\Res_{z_{0}}z_{0}^{-1}(z_{2}+z_{0})^{{\wt}a+n}z_2^{\wt b-n}
\langle u',Y(Y(a,z_{0})b,z_{2})u\rangle\\
& &=\<u',\sum_{k\in\Z_+}{\wt a+n\choose k}
(a_{k-1}b)(\wt(a_{k-1}b)-1)u\>\\
& &=\<u',\sum_{k\in\Z_+}{\wt a+n\choose k}o(a_{k-1}b)u\>\\
& &=\<u',o(\Res_zY(a,z)b\frac{(1+z)^{\wt a+n}}{z})u\>.
\end{eqnarray*}

On the other hand,  note that $b(\wt b-1+p)u=0$ if $p>n.$ So
\begin{eqnarray}
& &\ \ \ \Res_{z_{0}}z_{0}^{-1}(z_{0}+z_{2})^{{\wt}a+n}z_2^{\wt b-n}
\langle u',Y(a,z_{0}+z_{2})Y(b,z_{2})u\rangle\nonumber\\
& &=\<u',\sum_{i\in\Z_+}a(\wt a-2-i+n)\sum_{m\geq -n}b(\wt b-1-m)z_2^{-n+i+m}u\>\\
& &=\<u',\sum_{i=0}^na(\wt a-1-i)b(\wt b-1+i)u\>\\
& &+\<u',\sum_{i=1}^{n}a(\wt a-1+i)b(\wt b-1-i)u\>.\label{g6.12}
\end{eqnarray}

Note that the $A_n(V)$-module structure on $U$ is equivalent to 
\begin{eqnarray*}
& &o(a)o(b)u=a(\wt a-1)b(\wt b-1)u\\
& &=\sum_{m=0}^n(-1)^m{m+n\choose n}o(\Res_{z}Y(a,z)b
\frac{(1+z)^{\wt a+n}}{z^{m+n+1}})u.
\end{eqnarray*}
By (\ref{def}) with $s=2,$ $a_1=a,a_2=b,$ $p_1=k=-p_2$ ($k>0$) we see that
\begin{eqnarray}\label{a4.5}
& &\<u',o_k(a)o_{-k}(b)u\>=\<u',a(\wt a-1-k)b(\wt b-1+k)u\>\nonumber\\
& &=\<u',\sum_{m=0}^{n-k}(-1)^m{m+n+k\choose m}o(\Res_zY(a,z)b
\frac{(1+z)^{\wt a+n}}{z^{m+n+1+k}})u\>.
\end{eqnarray}
Thus 
\begin{eqnarray*}
& &\<u',\sum_{k=0}^na(\wt a-1-k)b(\wt b-1+k)u\>\\
& &=\<u',\sum_{k=0}^n\sum_{m=0}^{n-k}(-1)^m{m+n+k\choose m}o(\Res_{z}Y(a,z)b
\frac{(1+z)^{\wt a+n}}{z^{m+n+1+k}})u\>.
\end{eqnarray*}

Use Lie algebra bracket (\ref{ld}) to get
\begin{eqnarray*}
& &a(\wt a-1+k)b(\wt b-1-k)=b(\wt b-1-k)a(\wt a-1+k)\\
& &\ \ \ +\sum_{i\geq 0}{\wt a-1+k\choose i}
(a(i)b)(\wt a+\wt b-2-i).
\end{eqnarray*}
By (\ref{a4.5}),
\begin{eqnarray*}
& &\<u',b(\wt b-1-k)a(\wt a-1+k)u\>\\
& &=\sum_{m=0}^{n-k}(-1)^m{m+n+k\choose m}o(\Res_{z}Y(b,z)a
\frac{(1+z)^{\wt b+n}}{z^{m+n+1+k}})u\>.
\end{eqnarray*}
A proof similar to that of Lemma \ref{l2.2} (ii) shows that 
\begin{eqnarray*}
& &\sum_{m=0}^{n-k}(-1)^m{m+n+k\choose m}\Res_{z}Y(b,z)a\frac{(1+z)^{\wt b+n}}{z^{m+n+1+k}}\\
& &-\sum_{m=0}^{n-k}{m+n+k\choose m}(-1)^{n+k}\Res_zY(a,z)b\frac{(1+z)^{\wt a+m+k-1}}{z^{1+m+n+k}}\in O_n(V).
\end{eqnarray*}

We now have
\begin{eqnarray*}
& &\<u',\sum_{k=1}^{n}a(\wt a-1+k)b(\wt b-1-k)u\>\\
& &=\sum_{k=1}^{n}\sum_{m=0}^{n-k}{m+n+k\choose m}(-1)^{n+k}
\<u',o(\Res_zY(a,z)b\frac{(1+z)^{\wt a+m+k-1}}{z^{1+m+n+k}})u\>\\
& &\ \ \ +\sum_{k=1}^n\sum_{i\geq 0}{\wt a-1+k\choose i}
\<u',(a(i)b)(\wt a+\wt b-2-i)u\>\\
& &=\sum_{k=1}^{n}\sum_{m=0}^{n-k}{m+n+k\choose m}(-1)^{n+k}
\<u',o(\Res_zY(a,z)b\frac{(1+z)^{\wt a+m+k-1}}{z^{1+m+n+k}})u\>\\
& &\ \ \ +\sum_{k=1}^n\<u',o(\Res_zY(a,z)b(1+z)^{\wt a-1+k})u\>.
\end{eqnarray*}

So it is enough to show the following identity:
\begin{eqnarray*}
& &\sum_{k=0}^n\sum_{m=0}^{n-k}(-1)^m{m+n+k\choose m}
\frac{(1+z)^{\wt a+n}}{z^{m+n+1+k}}\\
& &+\sum_{k=1}^{n}\sum_{m=0}^{n-k}{m+n+k\choose m}(-1)^{n+k}\frac{(1+z)^{\wt a+m+k-1}}{z^{1+m+n+k}}\\
& &+\sum_{k=1}^n(1+z)^{\wt a-1+k}\\
& &=\frac{(1+z)^{\wt a+n}}{z},
\end{eqnarray*}
or equivalently,
\begin{eqnarray*}
& &\sum_{k=0}^n\sum_{m=0}^{n-k}(-1)^m{m+n+k\choose m}
\frac{(1+z)^n}{z^{m+n+k}}\\
& &+\sum_{k=1}^{n}\sum_{m=0}^{n-k}{m+n+k\choose m}(-1)^{n+k}\frac{(1+z)^{m+k-1}}{z^{m+n+k}}\\
& &=1.
\end{eqnarray*}
This identity is proved in Proposition \ref{ap1} in the Appendix. \qed

Proposition \ref{p3.3} is a consequence of the next lemma.
\bl{l6.8} For all $m\in\Z$ we have
\begin{eqnarray*}
& &\Res_{z_{0}}z_{0}^{m}(z_{0}+z_{2})^{{\wt}a+m+j}
\langle u',Y(a,z_{0}+z_{2})Y(b,z_{2})u\rangle\nonumber\\
&=&\Res_{z_{0}}z_{0}^{m}(z_{2}+z_{0})^{{\wt}a+m+j}
\langle u',Y(Y(a,z_{0})b,z_{2})u\rangle.
\end{eqnarray*}
\el

\pf This is true for $m\geq -1$ by Lemma \ref{l6.7}. Let us write $m=-k+i$ with
$i\in \Z_+$ and proceed by induction $k.$ Induction yields
\begin{eqnarray*}
& &\Res_{z_{0}}z_{0}^{-k}(z_{0}+z_{2})^{{\wt}a+m+j}
\langle u',Y(L(-1)a,z_{0}+z_{2})Y(b,z_{2})u\rangle\nonumber\\
&=&\Res_{z_{0}}z_{0}^{-k}(z_{2}+z_{0})^{{\wt}a+m+j}
\langle u',Y(Y(L(-1)a,z_{0})b,z_{2})u\rangle.\label{g6.14}
\end{eqnarray*}

Using the residue property $\Res_{z}f'(z)g(z)\!+\!\Res_zf(z)g'(z)\!=\!0$
and the $L(-1)$-derivation property $Y(L(-1)a,z)=\frac{d}{dz}Y(a,z)$ we have
\begin{eqnarray*}
& &\ \ \ \Res_{z_{0}}z_{0}^{-k}(z_{0}+z_{2})^{{\wt}a+1+m+j}
\langle u',Y(L(-1)a,z_{0}+z_{2})Y(b,z_{2})u\rangle\nonumber\\
& &=-\Res_{z_{0}}\left({\partial\over\partial z_{0}}
z_{0}^{-k}(z_{0}+z_{2})^{{\wt}a+1+m+j}\right)
\langle u',Y(a,z_{0}+z_{2})Y(b,z_{2})u\rangle\nonumber\\
& &=\Res_{z_{0}}k
z_{0}^{-k-1}(z_{0}+z_{2})^{{\wt}a+1+m+j}
\langle u',Y(a,z_{0}+z_{2})Y(b,z_{2})u\rangle\nonumber\\
& &\ \ \ -\Res_{z_{0}}({\wt}a+1+m+j)z_{0}^{-k}(z_{0}+z_{2})^{{\wt}a+
m+j}
\langle u',Y(a,z_{0}+z_{2})Y(b,z_{2})u\rangle\nonumber\\
& &=\Res_{z_{0}}k
z_{0}^{-k-1}z_{2}(z_{0}+z_{2})^{{\wt}a+m+j}
\langle u',Y(a,z_{0}+z_{2})Y(b,z_{2})u\rangle\nonumber\\
& &\ \ \ +\Res_{z_{0}}k
z_{0}^{-k}(z_{0}+z_{2})^{{\wt}a+m+j}
\langle u',Y(a,z_{0}+z_{2})Y(b,z_{2})u\rangle\nonumber\\
& &\ \ \ -\Res_{z_{0}}({\wt}a+1+m+j)z_{0}^{-k}(z_{2}+z_{0})^{{\wt}a+m+j}
\langle u',Y(Y(a,z_{0})b,z_{2})u\rangle\nonumber\\
& &=\Res_{z_{0}}k
z_{0}^{-k-1}z_{2}(z_{0}+z_{2})^{{\wt}a+m+j}
\langle u',Y(a,z_{0}+z_{2})Y(b,z_{2})u\rangle\nonumber\\
& &\ \ \ +\Res_{z_{0}}k
z_{0}^{-k}(z_{2}+z_{0})^{{\wt}a+m+j}
\langle u',Y(Y(a,z_{0})b,z_{2})u\rangle\nonumber\\
& &\ \ \ -\Res_{z_{0}}({\wt}a+1+m+j)z_{0}^{-k}(z_{2}+z_{0})^{{\wt}a+m+j}
\langle u',Y(Y(a,z_{0})b,z_{2})u\rangle,\label{g6.15}
\end{eqnarray*}
and
\begin{eqnarray*}
& &\Res_{z_{0}}z_{0}^{-k}(z_{2}+z_{0})^{{\wt}a+1+m+j}
\langle u',Y(Y(L(-1)a,z_{0})b,z_{2})u\rangle\nonumber\\
&=&-\Res_{z_{0}}\left({\partial\over\partial z_{0}}z_{0}^{-k}
(z_{2}+z_{0})^{{\wt}a+1+m+j}\right)\langle u',Y(Y(a,z_{0})b,z_{2})u\rangle
\nonumber\\
&=&\Res_{z_{0}}kz_{0}^{-k-1}(z_{2}+z_{0})^{{\wt}a+1+m+j}
\langle u',Y(Y(a,z_{0})b,z_{2})u\rangle\nonumber\\
& &-\Res_{z_{0}}({\wt}a+1+m+j)z_{0}^{-k}(z_{2}+z_{0})^{{\wt}a+m+j}
\langle u',Y(Y(a,z_{0})b,z_{2})u\rangle\nonumber\\
&=&\Res_{z_{0}}kz_{2}z_{0}^{-k-1}(z_{2}+z_{0})^{{\wt}a+m+j}
\langle u',Y(Y(a,z_{0})b,z_{2})u\rangle\nonumber\\
& &+\Res_{z_{0}}kz_{0}^{-k}(z_{2}+z_{0})^{{\wt}a+m+j}
\langle u',Y(Y(a,z_{0})b,z_{2})u\rangle\nonumber\\
& &-\Res_{z_{0}}({\wt}a+1+m+j)z_{0}^{-k}(z_{2}+z_{0})^{{\wt}a+m+j}
\langle u',Y(Y(a,z_{0})b,z_{2})u\rangle.
\end{eqnarray*}
This yields the identity:
\begin{eqnarray*}
& &\Res_{z_{0}}z_{0}^{-k-1}(z_{0}+z_{2})^{{\wt}a+m+j}
\langle u',Y(a,z_{0}+z_{2})Y(b,z_{2})u\rangle\nonumber\\
&=&\Res_{z_{0}}z_{0}^{-k-1}(z_{2}+z_{0})^{{\wt}a+m+j}
\langle u',Y(Y(a,z_{0})b,z_{2})u\rangle,\label{3.15}
\end{eqnarray*}
and the lemma is proved.\qed

Let us now introduce an arbitrary $\Z$-graded $\hat V$-module
$M=\oplus_{m\in \Z}M(m).$ As before we extend $M(n)^*$ to $M$ by letting it annihilate $M(m)$ for $m\ne n.$ The proof of Proposition of 6.1 
of [DLM] with $\<u',\cdot\>$ suitably inserted gives:

\begin{prop}\label{p3.4} Let $U$ be a subspace of $M(n)$ and $U'$
a subspace of $M(n)'$ such that

(i) $M=U(\hat V)U.$

(ii) For $a\in V$ and  $u\in U$
there is $k\in\Z$ such that 
\begin{eqnarray}\label{ea4}
\langle u',(z_{0}+z_{2})^{k+n}Y(a,z_{0}+z_{2})Y(b,z_{2})u\rangle
=\langle u',(z_{2}+z_{0})^{k+n}Y(Y(a,z_{0})b,z_{2})u\rangle  
\end{eqnarray}
for any $b\in V, u'\in U'$. Then in fact (\ref{ea4}) holds for any
$u\in M.$
\end{prop}

\begin{prop}\label{p3.5}
Let $M$ be as in Proposition \ref{p3.4}. Then for any $x\in
U(\hat V), a\in V, u\in M$, there is an integer $k$ such that 
\begin{eqnarray}\label{ea5}
\langle u',(z_{0}\!+\!z_{2})^{k+n}x\cdot Y(a,z_{0}\!+\!z_{2})Y(b,z_{2})u\rangle\!=\!\langle u',(z_{2}\!+\!z_{0})^{k+n}x\cdot Y(Y(a,z_{0})b,z_{2})u\rangle  
\end{eqnarray}
{\it for any $b\in V, u'\in U'$.}
\end{prop}

{\bf Proof.} Let $L$ be the subspace of $U(\hat V)$ consisting of those
$x$ for which  (\ref{ea5}) holds. Let $x\in L$, let $c$ be any homogeneous
element of $V,$ and let $m\in {\Z}$. Then from (\ref{3.6''}) we have
\begin{eqnarray}
& &\ \ \ \langle u',xc(m)Y(a,z_{0}+z_{2})Y(b,z_{2})u\rangle (z_{0}+z_{2})^{k+n}
\nonumber\\
& &=\sum_{i=0}^{\infty}{m\choose i}
(z_{0}+z_{2})^{k+n+n-i}\langle u',xY(c(i)a,z_{0}+z_{2})Y(b,z_{2})u\rangle
\nonumber\\
& &\ \ \ \ +\sum_{i=0}^{\infty}{m\choose i}
z_{2}^{n-i}(z_{0}+z_{2})^{k+n}\langle u',xY(a,z_{0}+z_{2})Y(c(i)b,z_{2})u\rangle\nonumber\\
& &\ \ \ \ +(z_{0}+z_{2})^{k+n}\langle u',xY(a,z_{0}+z_{2})Y(b,z_{2})c(m)u\rangle.
\end{eqnarray}
The same method that was used in the proof of Proposition \ref{p3.4} shows 
that  
$xc(m)\in L$. Since $U(\hat V)$ is generated by all such
$c(n)$'s, and since 
(\ref{ea5}) holds for $x=1$ by Proposition \ref{p3.4},
we conclude that $L=U(\hat V),$ as desired. \qed

We can now finish the proof of  Theorems \ref{t6.1} and Theorem \ref{t6.3}.
We can take $M=\bar M_n(U)$ in Proposition
\ref{p3.5}, as we may since $\bar M_n(U)$ certainly satisfies the conditions placed
on $M$ prior to Proposition \ref{p3.4} and in Proposition \ref{p3.4}. Then
{}from the definition of $W$ (\ref{g6.3}) 
and Propositions \ref{p3.3}, \ref{p3.4} and \ref{p3.5} we conclude
that $U(\hat V)W\subset J.$ It is clear that $L(U)$ is a quotient of
$\bar M_n(U)$  and hence an admissible $V$-module. Note that
$J\cap U=0.$ So $L(U)(n)$ contains $U$ as an
$A_n(V)$-submodule. 
This shows that  $\bar M_n(U)(n)\cong U$ as $A_n(V)$-modules.
If $\bar M_n(U)(0)=0$ then $U$ will be an $A_{n-1}(V)$-module, contradicting
the assumption on $U.$ This finishes the proof of Theorem \ref{t6.1}.
Theorem \ref{t6.3} is now obvious. 
\qed

At this point we have a pair of functors $\O_n,L_n$ defined on appropriate
module categories. It is clear that $\O_n/\O_{n-1}\circ L_n$ is equivalent to 
the identity. 

\bl{l7.1} Suppose that $U$ is a simple $A_n(V)$-module. 
Then $L_n(U)$ is a simple admissible $V$-module.
\el

\pf If $0\ne W\subset L_n(U)$ is an admissible 
submodule then, by the definition of $L_n(U)$, we have
$W(n)=W\cap L_n(U)(n)\ne 0$. As $W(n)$ is an 
$A_n(V)$-submodule of $U=L_n(U)(n)$ by Theorem \ref{t5.3} then
$U=W(n),$ whence $W\supset U(\hat V)W(n)=U(\hat V)U=L_n(U).$ \qed

\bt{t7.2} $L_n$ and $\O_n/O_{n-1}$
are equivalences when restricted to the full 
subcategories of completely reducible $A_n(V)$-modules whose irreducible
components cannot factor through $A_{n-1}(V)$ 
and 
completely reducible admissible $V$-modules respectively.
In particular, $L_n$ and $\Omega_n/\O_{n-1}$ induce mutually
inverse bijections on the isomorphism classes of simple objects
in the category of $A_n(V)$-modules which cannot factor through $A_{n-1}(V)$
and admissible $V$-modules
respectively.
\end{thm}

\pf We have $\Omega_n/O_{n-1}(L(U))\cong U$ for any $A_n(V)$-module by Theorem
\ref{t6.3}.

If $M$ is a completely reducible admissible 
$V$-module we must show $L_n(\O_n/\O_{n-1}(M))\cong M.$ For this we may take
$M$ simple, whence $\O_n/\O_{n-1}(M)$ is simple by Proposition 
\ref{l2.9} (ii) and then $L_n(\O_n/\O_{n-1}(M))$ is simple by Lemma \ref{l7.1}. Since both
$M$ and $L_n(\Omega_n/\O_{n-1}(M))$ are simple quotients
of the universal object $\bar{M}_n(\Omega_n/\O_{n-1}(M))$ then they are isomorphic
by Theorems \ref{t6.1} and \ref{t6.3}. \qed

The following theorem is a generalization of Theorem 8.1 of [DLM].

\bt{t8.1} Suppose that $V$ is a rational vertex operator algebra. Then the 
following hold:

(a) $A_n(V)$ is a finite-dimensional, semisimple associative algebra.

(b) The functors $L_n,\O_n/O_{n-1}$ are mutually inverse categorical equivalences
between the category of $A_n(V)$-modules whose irreducible components
cannot factor through $A_{n-1}(V)$ and the category of admissible 
$V$-modules.

(c) The functors $L_n,\O_n/\O_{n-1}$ induce mutually inverse categorical equivalences
between the category of finite-dimensional 
$A_n(V)$-modules whose irreducible components
cannot factor through $A_{n-1}(V)$ and the category of ordinary $V$-modules.
\et

\pf (b) follows from Theorem \ref{t7.2} and (a).
Since $V$ is rational any irreducible admissible $V$-module
is an ordinary module by Theorem 8.1 of [DLM]. Now (c) follows from (b).
It remains to prove (i). 

Let $W$ be an $A_n(V)$-module. Set $U=W\oplus V(n).$ Then $U$ is 
an $A_n(V)$-module which cannot factor through $A_{n-1}(V).$ Now $L_n(U)$
is admissible and hence a direct sum of irreducible ordinary $V$-modules.
Thus $\O_n(L_n(U))/\O_{n-1}(L_n(U))\simeq U$ is a direct sum of
finite-dimensional irreducible $A_n(V)$-modules and so is $W.$ 
\qed

It is believed that if $A(V)=A_0(V)$ is semisimple then $V$ is rational.
We cannot solve this problem completely in this paper. But we have some
partial results which are applications of $A_n(V)$-theory.

\begin{thm}\label{t4.13} If all $A_n(V)$ are finite-dimensional 
semisimple algebras then $V$ is rational.
\end{thm}

\pf Since $A(V)$ is semisimple $V$ has only finitely many irreducible
admissible modules which are necessarily ordinary $V$-modules.
For any $\lambda\in\C$ let ${\cal M}_{\l}$ be the set of irreducible
admissible modules whose weights are congruent to $\l$ module $\Z.$
Then for each $W\in {\cal W}_{\l}$ we have 
$W=\oplus_{n\in\Z_+}W_{\l+n_W+n}=\oplus_{n\in\Z_+}W(n)$
where  $n_W\in\Z$ and $W_{\l+n_W+n}=W(n).$ 
Since $L(-1): W(n)\to W(n+1)$ is injective if $n$ is large (see [L1])
there exists an $m_{\l}\in \N$ such that the weight
space $W_{\l+m}\ne 0$ for any $W\in {\cal W}_{\l}$ and $m\geq m_{\l}.$

Consider any admissible module $M$ whose weights are
in $\l+Z$ and whose homogeneous subspace $M_{\l+m}$ with some $m\geq m_{\l}$
is 0. Let $U$ be an irreducible $A(V)$-submodule of $M(0).$ Then
$L_0(U)=L(U)$ is an irreducible $V$-module such that 
$L(U)(0)=U$ and $L(U)_{\l+m}=0.$ Thus $L(U)=0$ and $U=0.$ This implies
that $M=0.$ 

Now take an admissible module $M=\oplus_{k\in\Z_+}M(k).$ Then $M(0)$ is a
direct sum of simple $A(V)$-modules as $A(V)$ is semisimple.  Let $U$
be an $A(V)$-submodule of $M(0)$ isomorphic to $W(0)=W_{\l+n_W}$ for
some $W\in {\cal M}_{\l}.$ We assert that the submodule $N$ of $M$
generated by $U$ is irreducible and necessarily isomorphic to $W.$
First note that $N$ has an irreducible quotient isomorphic to $W.$
Take $n\in \N$ such that $n+n_W\geq m_{\l}.$ Observe that $\bar
M_n(W(n))/\bar J=L_n(W(n))$ is isomorphic to $W$ where $\bar J$ is a
maximal submodule of $\bar M_n(W(n))$ such that $\bar J\cap W(n)=0.$
Since $\bar J_{\l+n_W+n}=0$ we see that $\bar J=0$ and $\bar
M_n(W(n))=L_n(W(n))\simeq W.$ Write $N(n)$ as a direct sum of $W(n)$
and another $A_n(V)$-submodule $N(n)'$ of $N(n)$ as $A_n(V)$ is
semisimple. Clearly the submodule of $N(n)$ generated by
$W(n)$ is isomorphic to $W.$ This shows that $N$ must be isomorphic to
$W,$ as claimed.

It is obvious now that the submodule $U(\hat V)M(0)$ generated
by $M(0)$ is completely reducible. Use the semisimplicity of $A_1(V)$
we can decompose $M(1)$ into a direct sum of $A_1(V)$-modules
$(U(\hat V)M(0))(1)\oplus M(1)'.$ The same argument shows that $U(\hat V)M(1)'$
is a completely reducible submodule of $M.$ Continuing in this way proves
that $M$ is completely reducible. \qed

\begin{rem} From the proof of Theorem \ref{t4.13}, we see, in fact, that
we can weaken the assumption in Theorem \ref{t4.13}. Namely we only
need to assume that $A_n(V)$ is semisimple if $n$ is large. 
\end{rem}

\section{Appendix}

In this appendix  we prove several combinatorial identities which are used 
in the previous sections. 

For $n\geq 0$ define 
\begin{eqnarray*}
& &A_n(z)=\sum_{k=0}^n\sum_{m=0}^{n-k}(-1)^m{m+n+k\choose m}
\frac{(1+z)^n}{z^{m+n+k}}\\
& &\ \ \ \ \ +\sum_{k=1}^{n}\sum_{m=0}^{n-k}{m+n+k\choose m}(-1)^{n+k}\frac{(1+z)^{m+k-1}}{z^{m+n+k}}.
\end{eqnarray*}
Using the well-known identity 
$$\sum_{k=0}^{i}(-1)^{k}{n\choose k}
=(-1)^{i}{n-1\choose i}$$ 
we can rewrite $A_n(z)$ as
\begin{eqnarray*}
& &A_n(z)=\sum_{k=0}^n\sum_{m=0}^{k}(-1)^m{n+k\choose m}
\frac{(1+z)^n}{z^{n+k}}+\sum_{k=1}^{n}\sum_{m=0}^{k-1}{n+k\choose m}(-1)^{n+k+m}\frac{(1+z)^{k-1}}{z^{n+k}}\\
& &=\sum_{k=0}^n(-1)^k{n+k-1\choose k}\frac{(1+z)^n}{z^{n+k}}
-(-1)^n\sum_{k=1}^{n}{n+k-1\choose k-1}\frac{(1+z)^{k-1}}{z^{n+k}}.
\end{eqnarray*}
\begin{prop}\label{ap1} $A_n(z)=1$ for all $n\geq 0.$
\end{prop}

\pf Set
$$B_n(z)=\sum_{k=0}^n(-1)^k{n+k-1\choose k}\frac{(1+z)^n}{z^{n+k}}$$
$$C_n(z)=\sum_{k=1}^{n}{n+k-1\choose k-1}\frac{(1+z)^{k-1}}{z^{n+k}}.$$
Then
\begin{eqnarray*}
& &B_n(z)=\sum_{k=0}^{n-1}(-1)^k\left({n+k-2\choose k}
+{n+k-2\choose k-1}\right)\frac{(1+z)^n}{z^{n+k}}+(-1)^n{2n-1\choose n}\frac{(1+z)^n}{z^{2n}}\\
& &\ \ \ \ \ =\frac{1+z}{z}B_{n-1}(z)
+\sum_{k=0}^{n-2}(-1)^{k+1}{n+k-1\choose k}
\frac{(1+z)^n}{z^{n+k+1}}+(-1)^n{2n-1\choose n}\frac{(1+z)^n}{z^{2n}}\\
& &\ \ \ \ \ =\frac{1+z}{z}B_{n-1}(z)-\frac{1}{z}B_{n}(z)+(-1)^{n-1}{2n-2\choose n-1}\frac{(1+z)^n}{z^{2n}}\\
& &\ \ \ \ \ \ \ \ \ \ +(-1)^{n}{2n-1\choose n}\frac{(1+z)^{n+1}}{z^{2n+1}}.
\end{eqnarray*}
Solving $B_n(z)$ gives 
$$B_n(z)=B_{n-1}(z)+(-1)^{n-1}\frac{(1+z)^{n-1}}{z^{2n-1}}
{2(n-1)\choose n-1}+(-1)^n\frac{(1+z)^{n}}{z^{2n}}{2n-1\choose n}.$$
Similarly,
\begin{eqnarray*}
& &C_n(z)=\sum_{k=1}^{n-1}\left({n+k-2\choose k-1}
+{n+k-2\choose k-2}\right)\frac{(1+z)^{k-1}}{z^{n+k}}+{2n-1\choose n-1}\frac{(1+z)^{n-1}}{z^{2n}}\\
& &\ \ \ \ \ =\frac{1}{z}C_{n-1}(z)
+\sum_{k=0}^{n-2}{n+k-1\choose k-1}
\frac{(1+z)^k}{z^{n+k+1}}+{2n-1\choose n-1}\frac{(1+z)^{n-1}}{z^{2n}}\\
& &\ \ \ \ \ =\frac{1}{z}C_{n-1}(z)+\frac{1+z}{z}C_{n}(z)+{2n-2\choose n-1}\frac{(1+z)^{n-1}}{z^{2n}}\\
& &\ \ \ \ \ \ \ \ \ \ -{2n-1\choose n-1}\frac{(1+z)^{n}}{z^{2n+1}}.
\end{eqnarray*}
Thus
\begin{eqnarray*}
& &(-1)^{n+1}C_n(z)=(-1)^nC_{n-1}(z)+(-1)^n{2n-2\choose n-1}
\frac{(1+z)^{n-1}}{z^{2n-1}}\\
& &\ \ \ \ \ \ \ \ \ \ +(-1)^{n-1}{2n-1\choose n-1}\frac{(1+z)^{n}}{z^{2n}}.
\end{eqnarray*}
Thus
\begin{eqnarray*}
A_n(z)=B_n(z)+(-1)^{n-1}C_n(z)=A_{n-1}.
\end{eqnarray*}
Note that $A_0(z)=1$ and the proposition follows.
\qed 

For $n\geq 0$ we define 
$$F_n(z)=\sum_{m=0}^n{m+n\choose
n}\frac{(-1)^m(1+z)^{n+1}-(-1)^n(1+z)^m}{z^{n+m+1}}.$$ 
\begin{prop}\label{ap2} $F_n(z)=1$ for all $n.$
\end{prop}
\pf Set 
$$D_n(z)=\sum_{m=0}^n{m+n\choose
n}(-1)^m\frac{(1+z)^{n+1}}{z^{n+m+1}}$$
$$E_n(z)=\sum_{m=0}^n{m+n\choose
n}(-1)^m\frac{(1+z)^m}{z^{n+m+1}}.$$
Then 
\begin{eqnarray*}
& &D_n(z)=B_{n+1}(z)+(-1)^n{2n+1\choose n}\frac{(1+z)^{n+1}}{z^{2n+2}}\\
& &\ \ \ =B_n(z)+(-1)^{n}\frac{(1+z)^{n}}{z^{2n+1}}
{2n\choose n}+(-1)^{n+1}\frac{(1+z)^{n+1}}{z^{2n+2}}{2n+1\choose n+1}\\
& &\ \ \ \ \ \ +(-1)^n{2n+1\choose n}\frac{(1+z)^{n+1}}{z^{2n+2}}\\
& &\ \ \ =B_n(z)+(-1)^{n}{2n\choose n}\frac{(1+z)^{n}}{z^{2n+1}}
\end{eqnarray*}
and
\begin{eqnarray*}
E_n(z)=C_{n}(z)+{2n\choose n}\frac{(1+z)^{n}}{z^{2n+1}}.
\end{eqnarray*}
Thus
$$ F_n(z)=D_n(z)+(-1)^{n+1}E_n(z)=A_n(z)=1,$$
as required.  \qed

For $n\geq 0$ define 
$$a_n(w,z)=\sum_{m=0}^n(-1)^{m}{m+n\choose
n}\left(\sum_{i=0}^{n-m}\sum_{j\geq 0}
{-m-n-1\choose i}{m\choose j}(-1)^i\frac{w^{i+j}}{z^{i+m}}-1\right).$$
Note that if $p>0, k>0$ then ${-p\choose k}=(-1)^k{p+k-1\choose k}.$ We can 
rewrite $a_n(w,z)$ as 
$$a_n(w,z)=\sum_{m=0}^n(-1)^{m}{m+n\choose
n}\left(\sum_{i=0}^{n-m}\sum_{j\geq 0}
{m+n+i\choose i}{m\choose j}\frac{w^{i+j}}{z^{i+m}}-1\right).$$

\begin{prop}\label{ap3} The $a_n(w,z)=0$ for all $n\geq 0.$
\end{prop}

\pf Regarding $a_n(w,z)$ as  a polynomial in $z^{-1},$
the coefficient of $z^{-p}$ in $A_n(w,z)$ ($0\leq p\leq n$) is equal to
(setting $m+i=p$)
\begin{eqnarray*}
& &\sum_{m=0}^p(-1)^{m}{m+n\choose
n}\sum_{j\geq 0}
{n+p\choose p-m}{m\choose j}w^{p-m+j}-(-1)^p{p+n\choose n}\\
& &\ \ \ 
=w^p\sum_{m=0}^p(-1)^{m}{m+n\choose
n}{n+p\choose p-m}(1+1/w)^m-(-1)^p{p+n\choose n}.
\end{eqnarray*}
So the coefficient of $z^{-p}w^0$ in $a_n(w,z)$ equals 0. 

If $0\leq q\leq p-1,$ the coefficient of $z^{-p}w^{p-q}$ in $a_n(w,z)$
is equals to 
$$c_n(p,q)=\sum_{m=0}^p(-1)^{m}{m+n\choose
n}{n+p\choose n+m}{m\choose q}$$
which is defined for any $n,p,q\geq 0.$ 
So we must prove that $a_n(p,q)=0$ for $1\leq q+1\leq p\leq n.$
Recall ${l\choose k}={l-1\choose k}+{l-1\choose k-1}.$ Then $c_n(p,q)$ is 
equal to 
\begin{eqnarray*}
& &\ \ \   \sum_{m=0}^p(-1)^{m}{m+n\choose
n}\left({n+p-1\choose n+m}+{n+p-1\choose n+m-1}\right){m\choose q}\\
& &=(-1)^p{p+n\choose n}{p\choose q}+\sum_{m=0}^{p-1}(-1)^{m}{m+n\choose
n}{n+p-1\choose n+m-1}{m\choose q}+c_n(p-1,q)\\
& &=\sum_{m=0}^{p-1}(-1)^{m}\left({m+n-1\choose
n-1}+{m+n-1\choose
n}\right){n+p-1\choose n+m-1}{m\choose q}\\
& &\ \ \ \ \ +(-1)^p{p+n\choose n}{p\choose q}+c_n(p-1,q)\\
& &=c_{n-1}(p,q)+c_n(p-1,q)-(-1)^p{p+n-1\choose n-1}{p\choose q}+(-1)^p{p+n\choose n}{p\choose q}\\
& &\ \ \ \ -\sum_{m=0}^{p-1}(-1)^{m-1}{m+n-1\choose
n}{n+p-1\choose n+m-1}\left({m-1\choose q}+{m-1\choose q-1}\right)\\
& &=c_{n-1}(p,q)+c_n(p-1,q)+(-1)^p{p+n-1\choose n}{p\choose q}\\
& &\ \ \ \ -\sum_{m-1=0}^{p-2}(-1)^{m-1}{m+n-1\choose
n}{n+p-1\choose n+m-1}\left({m-1\choose q}+{m-1\choose q-1}\right)\\
& &=c_{n-1}(p,q)-c_n(p-1,q-1)+(-1)^p{p+n-1\choose n}{p\choose q}\\
& &\ \ \ \ \ \ 
+(-1)^{p-1}{p-1+n\choose n}\left({p-1\choose q}+{p-1\choose q-1}\right)
\\
& &=c_{n-1}(p,q)-c_n(p-1,q-1).
\end{eqnarray*}
That is,
$$ c_n(p,q)=c_{n-1}(p,q)-c_n(p-1,q-1).$$
So by induction it is enough to show that
$c_0(p,q)=0$ and $c_n(p,0)=0$ if $p>q.$ 
But this is clear from the definition.
\qed

\end{document}